\documentclass[a4paper,11pt,reqno]{amsart}
\pdfoutput=1
\usepackage{amssymb,amsmath,amsthm,amsfonts,centernot}
\usepackage[colorlinks=true,linkcolor=blue,citecolor=blue,urlcolor=cyan, pdfpagelabels=false]{hyperref}
\usepackage{amsmath,tikz}
\usepackage{makecell}
\usepackage{dsfont}
\usetikzlibrary{matrix}
\usepackage{perpage} 
\usepackage{makecell}
\MakePerPage{footnote}
\begin{document}
\numberwithin{equation}{section}
\newtheorem{theorem}{Theorem}
\newtheorem{algo}{Algorithm}
\newtheorem{lem}{Lemma} 
\newtheorem{de}{Definition} 
\newtheorem{ex}{Example}
\newtheorem{pr}{Proposition} 
\newtheorem{claim}{Claim} 
\newtheorem{re}{Remark}
\newtheorem{co}{Corollary}
\newtheorem{conv}{Convention}
\newcommand{\di}{\hspace{1.5pt} \big|\hspace{1.5pt}}
\newcommand{\idi}{\hspace{.5pt}|\hspace{.5pt}}
\newcommand{\hs}{\hspace{1.3pt}}
\newcommand{\thmf}{Theorem~1.15$'$}
\newcommand{\ndi}{\centernot{\big|}}
\newcommand{\nidi}{\hspace{.5pt}\centernot{|}\hspace{.5pt}}
\newcommand{\lp}{\mbox{$\hspace{0.12em}\shortmid\hspace{-0.62em}\alpha$}}
\newcommand{\PQ}{\bb{P}^1(\bb{Q})}
\newcommand{\pmn}{\cl{P}_{M,N}}
\newcommand{\lcm}{\operatorname{lcm}}
\newcommand{\he}{holomorphic eta quotient\hspace*{2.5pt}}
\newcommand{\hes}{holomorphic eta quotients\hspace*{2.5pt}}
\newcommand{\defG}{Let $G\subset\GG$ be a subgroup that is conjugate to a finite index subgroup of $\G$. } 
\newcommand{\defg}{Let $G\subset\GG$ be a subgroup that is conjugate to a finite index subgroup of $\G$\hs\hs} 
\renewcommand{\phi}{\varphi}
\newcommand{\Z}{\bb{Z}}
\newcommand{\ZD}{\Z^{\D}}
\newcommand{\N}{\bb{N}}
\newcommand{\Q}{\bb{Q}}
\newcommand{\pii}{{{\pi}}}
\newcommand{\R}{\bb{R}}
\newcommand{\C}{\bb{C}}
\newcommand{\I}{\hs\cl{I}_{n,N}}
\newcommand{\SL}{\operatorname{SL}_2(\Z)}
\newcommand{\St}{\operatorname{Stab}}
\newcommand{\D}{\cl{D}_N}
\newcommand{\rh}{{{\boldsymbol\rho}}}
\newcommand{\bh}{{\cl{M}}} 
\newcommand{\lv}{\hyperlink{level}{{\text{level}}}\hspace*{2.5pt}}
\newcommand{\fct}{\hyperlink{factor}{{\text{factor}}}\hspace*{2.5pt}}
\newcommand{\q}{\hyperlink{q}{{\mathbin{q}}}}
\newcommand{\rd}{\hyperlink{redu}{{{\text{reducible}}}}\hspace*{2.5pt}}
\newcommand{\ird}{\hyperlink{irredu}{{{\text{irreducible}}}}\hspace*{2.5pt}}
\newcommand{\str}{\hyperlink{strong}{{{\text{strongly reducible}}}}\hspace*{2.5pt}}
\newcommand{\rdn}{\hyperlink{redon}{{{\text{reducible on}}}}\hspace*{2.5pt}}
\newcommand{\atl}{\hyperlink{atinv}{{\text{Atkin-Lehner involution}}}\hspace*{3.5pt}}
\newcommand{\atls}{\hyperlink{atinv}{{\text{Atkin-Lehner involutions}}}\hspace*{3.5pt}}
\newcommand{\T}{\mathrm{T}}
\renewcommand{\H}{\fr{H}}
\newcommand{\W}{\text{\calligra W}_n}
\newcommand{\GG}{\cl{G}}
\newcommand{\g}{\fr{g}}
\newcommand{\Gm}{\Gamma}
\newcommand{\Gmtl}{\widetilde{\Gamma}_\ell}
\newcommand{\gm}{\gamma}
\newcommand{\go}{\gamma_1}
\newcommand{\gmt}{\widetilde{\gamma}}
\newcommand{\gmdt}{\widetilde{\gamma}'}
\newcommand{\gmot}{\widetilde{\gamma}_1}
\newcommand{\gmdot}{{\widetilde{\gamma}}'_1}
\newcommand{\s}{\Large\text{{\calligra r}}\hspace{1.5pt}}
\newcommand{\ms}{m_{{{S}}}}
\newcommand{\nisim}{\centernot{\sim}}
\newcommand{\level}{\hyperlink{level}{{\text{level}}}}
\newcommand{\Redcon}{the \hyperlink{red}{\text{Reducibility~Conjecture}}}
\newcommand{\ReDcon}{The \hyperlink{red}{\text{Reducibility~Conjecture}}}
\newtheorem*{ThmA}{Theorem A}
\newtheorem{Conj}{Conjecture}
\newtheorem*{ThmB}{Theorem B}
\newtheorem*{ThmC}{Theorem C}
\newtheorem*{lmB}{Lemma B}
\newtheorem*{CorA}{Corollary A}
\newtheorem*{CorB}{Corollary B}
\newtheorem*{CorC}{Corollary C}
\newcommand{\Conred}{Conjecture~$1'$}
\newcommand{\effth}{Theorem~$\ref{27.11.2015}'$}
\newcommand{\Conredd}{Conjecture~$1''$}
\newcommand{\Conreddd}{Conjecture~$1'''$}
\newcommand{\Conired}{Conjecture~$2'$}
\newtheorem*{pro}{\textnormal{\textit{Proof of Lemma~\ref{27.11.2015.1}}}}
\newtheorem*{cau}{Caution}
\newtheorem{thrmm}{Theorem}[section]
\newtheorem*{thmA}{Theorem~A}
\newtheorem*{corA}{Corollary}
\newtheorem*{corB}{Corollary~B}
\newtheorem*{corC}{Corollary~C}
\newtheorem{no}{Notation}
\renewcommand{\thefootnote}{\fnsymbol{footnote}}
\newtheorem{oq}{Question}
\newtheorem{hy}{Hypothesis} 
\newtheorem{expl}{Example}
\newcommand\ileg[2]{\bigl(\frac{#1}{#2}\bigr)}
\newcommand\leg[2]{\Bigl(\frac{#1}{#2}\Bigr)}
\newcommand{\e}{\eta}
\newcommand{\sgn}{\operatorname{sgn}}
\newcommand{\bb}{\mathbb}
\newtheorem*{conred}{Conjecture~\ref{con1}$\mathbf{'}$}
\newtheorem*{conredd}{Conjecture~\ref{con1}$\mathbf{''}$}
\newtheorem*{conreddd}{Conjecture~\ref{con1}$\mathbf{'''}$}
\newtheorem*{conired}{Conjecture~\ref{19.1Aug}$\mathbf{'}$}
\newtheorem*{efth}{Theorem~\ref{27.11.2015}$\mathbf{'}$}
\newtheorem*{eflem}{Theorem~\ref{15May15}.$(b)\mathbf{'}$}
\newtheorem*{procl}{\textnormal{\textit{Proof}}}
\newtheorem*{thmff}{Theorem~\ref{7.4Jul}$\mathbf{'}$}
\newtheorem*{coo}{Corollary~\ref{17Aug}$\mathbf{'}$}
\newcommand{\cooo}{Corollary~\ref{17Aug}$'$}
\newtheorem*{cotw}{Corollary~\ref{17.1Aug}$\mathbf{'}$}
\newcommand{\cotww}{Corollary~\ref{17.1Aug}$'$}
\newtheorem*{cothr}{Corollary~\ref{17.2Aug}$\mathbf{'}$}
\newcommand{\cothre}{Corollary~\ref{17.2Aug}$'$}
\newtheorem*{cne}{Corollary~\ref{15.5Aug}$\mathbf{'}$}
\newcommand{\cnew}{Corollary~\ref{15.5Aug}$'$\hspace{3.5pt}}
\newcommand{\fr}{\mathfrak}
\newcommand{\cl}{\mathcal}
\newcommand{\lr}[1]{\left(#1\right)}
\newcommand{\rad}{\mathrm{rad}}
\newcommand{\ord}{\operatorname{ord}}
\newcommand{\m}{\setminus}
\newcommand{\G}{\Gamma_1}
\newcommand{\GN}{\Gamma_0(N)}
\newcommand{\X}{\widetilde{X}}
\renewcommand{\P}{{\textup{p}}} 
\newcommand{\al}{{\hs\operatorname{al}}}
\newcommand{\p}{p_\text{\tiny (\textit{N})}}
\newcommand{\pN}{p_\text{\tiny\textit{N}}}
\newcommand{\bt}{\mbox{$\raisebox{-0.59ex}
  {${{l}}$}\hspace{-0.215em}\beta\hspace{-0.88em}\raisebox{-0.98ex}{\scalebox{2}
  {$\color{white}.$}}\hspace{-0.416em}\raisebox{+0.88ex}
  {$\color{white}.$}\hspace{0.46em}$}{}}
  \newcommand{\un}{\hs\underline{\hspace{5pt}}\hs}
\newcommand{\U}{u_\textit{\tiny N}}
\newcommand{\Upr}{u_{\textit{\tiny N}^\prime}}
\newcommand{\Up}{u_{\textit{\tiny p}^\textit{\tiny e}}}
\newcommand{\Un}{u_{\textit{\tiny p}_\textit{\tiny 1}^{\textit{\tiny e}_\textit{\tiny 1}}}}
\newcommand{\Um}{u_{\textit{\tiny p}_\textit{\tiny m}^{\textit{\tiny e}_\textit{\tiny m}}}}
\newcommand{\Ut}{u_{\text{\tiny 2}^\textit{\tiny a}}}
\newcommand{\At}{A_{\text{\tiny 2}^\textit{\tiny a}}}
\newcommand{\Uh}{u_{\text{\tiny 3}^\textit{\tiny b}}}
\newcommand{\Ah}{A_{\text{\tiny 3}^\textit{\tiny b}}}
\newcommand{\Uprl}{u_{\textit{\tiny N}_1}}
\newcommand{\Uprlm}{u_{\textit{\tiny N}_i}}
\newcommand{\UM}{u_\textit{\tiny M}}
\newcommand{\UMp}{u_{\textit{\tiny M}_1}}
\newcommand{\w}{\omega_\textit{\tiny N}}
\newcommand{\wm}{\omega_\textit{\tiny M}}
\newcommand{\wa}{\omega_{\text{\tiny N}_\textit{\tiny a}}}
\newcommand{\wma}{\omega_{\text{\tiny M}_\textit{\tiny a}}}
\renewcommand{\P}{{\textup{p}}}


\title[Primes and polygonal numbers]{Primes and polygonal numbers}

\author{Soumya Bhattacharya and Habibur Rahaman}
\address
{Department of Mathematics and Statistics, Indian Institute of Science Education and Research Kolkata, Mohanpur, West Bengal - 741246, India.} 
\email{soumya.bhattacharya@iiserkol.ac.in, hr21rs044@iiserkol.ac.in}


\maketitle
 \begin{abstract}

A linear combination $aT_r(m)+bT_s(n)$ 
of an \mbox{$r$-gonal} number $T_r(m)$ and an $s$-gonal number $T_s(n)$ with mutually coprime positive
integer coefficients $a$ and $b$ produces infinitely many primes as $m$ and~$n$ varies over the natural numbers, whereas 
the 
sum of the reciprocals of 
such primes 
converges unless
$T_r(m)=m^2$ and $T_s(n)=n^2$.
For each pair of 
coprime positive integers $a$ and $b$,
there are arbitrary long arithmetic progressions among the primes of the form
$am^2+bn^2$.
 \end{abstract}

\section{Introduction and the main result}

A prime is a sum of two squares if and only if it is of the form $4k+1$. Since there are infinitely many primes of the form $4k+1$, it follows that there are infinitely many primes which are sum of two squares. 
Since the sums of couples of triangular numbers are intimately related  
with the sums of couples of squares, 
viz. $n$ is a sum of two triangular numbers if and only if
$4n+1$ is a sum of two squares, it is natural to expect that the primes which are sums of two triangular numbers will also have an equally appealing description in terms of some arithmetic progressions. 
However, it turns out not to be the case. In fact, the sum of the reciprocals of such primes converges. 
One may still wonder whether there are infinitely many such primes
or more generally, 
whether for a given 
pair of 
integers $r,s>2$, there are infinitely many primes which 
are sum of an $r$-gonal number and an $s$-gonal number,
where the $m$-th $r$-gonal number is defined by (see~\cite{CG})
\begin{equation}
    T_r(m) \; := \; (r-2)\frac{m(m-1)}{2}+m.\label{defpolynum}
\end{equation}
Even more generally, one may ask whether for each pair of integers $r,s>2$ and for 
every pair of mutually coprime positive integers $a$ and $b$, there are infinitely many primes
of the form $a T_r(m)+b T_s(n).$
We answer this in the affirmative. 
\begin{theorem}\label{thm_2}
Given a pair of 
 integers $r,s>2$ and a pair of
 mutually coprime positive integers $a$ and $b$, for all $m,n\in\N$, the number of primes of the form \begin{align}\label{def_ uT+vT}
     a T_r(m)+b T_s(n)
 \end{align}
which are less than $N$ is 
$$\asymp\begin{cases}
\dfrac{N}{\log N}, &\text{ if } r=s=4\\ \\
\dfrac{N}{(\log N)^{3/2}} &\text{otherwise.}
\end{cases}$$
\end{theorem}
Recall that $f\asymp g$ if and only if there exists a pair of constants $C,C'>0$ such that $Cg\leq f \leq C'g$.
In the first case above, we have $T_4(m)=m^2$.
The surprise is 
in 
the other case, where
even if one of the summands is a 
multiple of a square,
the corresponding set of primes has a lower rate of growth.
From the above theorem, via Abel's summation formula, we obtain
\begin{co}\label{co0}
Given $a,b,r,s\in\N$ with $r,s>2$ and $\gcd(a,b)=1$,
the sum of the reciprocals of the primes of the form 
$a T_r(m)+b T_s(n)$
diverges if and only if $r=s=4$. 
\end{co}

In particular, 
Theorem~\ref{thm_2}
answers 
a question
from 
\cite{B}, which asks how many primes are sums of two triangular numbers. 

 \begin{co}
There are infinitely many primes which are sum of two triangular numbers.
\end{co}

In fact, there are infinitely many such primes in 
each of the 
residue classes $\pm1$ modulo $4$.

\begin{theorem}\label{thm_3}
Given two mutually coprime positive integers $a$ and $b$,
the number of the primes in each of the 
residue classes 
$\pm1$ modulo $4$,
which are less than $N$ 
and
which are 
linear combinations of
two odd-sided polygonal numbers 
with coefficients $a$ and $b$
is 
$$\asymp \dfrac{N}{(\log N)^{3/2}}.$$ 
\end{theorem}
In particular, for $a=b=1$ and for the primes congruent to $1$ modulo $4$, we have the following.
\begin{corA}\label{cor2_thm3}
    For each odd integer $r>1$, the number of primes $p\leq N$, which are sums of two $r$-gonal numbers and also sums of two squares, is 
       $$ \asymp \frac{N}{(\log N)^{3/2}}.$$
\end{corA}

Since a number $n$ is a sum of 
two triangular numbers
if and only if $4n+1$ is a
sum of two squares\footnote{
It follows from \eqref{defpolynum} that for all $r$, we have 
    \begin{align*}
        4(r-2)(T_r(m)+T_r(n))+(r-4)^2 \; = \; \lr{(r-2)(m-n)}^2+\lr{(r-2)(m+n)-r+4}^2.
    \end{align*}
    Conversely, 
    if $4(r-2)M\; + \;(r-4)^2 \; = \; A^2\; + \;B^2$ such that $A\; \pm \; B \; \equiv \; r\pmod{2(r-2)}$, then $M$ is a sum of two $r$-gonal numbers.
    } (see~\cite{B}) and since a prime is a sum of two squares if and only if it is of the form $4k+1$ (see~\cite{Z}),
it follows that 
\begin{corA}\label{cor1_thm3}
The number of primes $p\leq N$ such that 
both $p$ and $4p+1$ are sums of two squares is
$$\asymp \dfrac{N}{(\log N)^{3/2}}.$$ 
\end{corA}
By a \emph{primitive arithmetic progression}, we mean 
an integer arithmetic progression that has a
pair of mutually coprime terms. Now, looking at Theorem~\ref{thm_3},
one may wonder whether every primitive arithmetic progression contains 
infinitely many primes which are sums of two triangular numbers.
However, that does not hold in general!
For example, if $n$ is congruent to $5$ or $8$ modulo 9, then
$4n+1$ is not a sum of two squares. 
So, none of
the arithmetic progressions $\{9k+5\}_{k\in \N}$ and $\{9k+8\}_{k\in \N}$ 
contains 
a sum of two triangular numbers. 
This prompts us to 
refine our query as follows:
\begin{oq}Given $a,b,r,s\in\N$ with $r,s>2$ and $\gcd(a,b)=1$,
if there is no 
congruence obstruction, does every primitive arithmetic progress\-ion contain infinitely many primes which are 
of the form $aT_r(m)+bT_s(n)?$
\end{oq}
\vspace*{2pt}

Theorem~\ref{thm_2} together with 
the Prime Number Theorem for arithmetic progressions (Corollary 11.20 in \cite{MV})
imply the following. 
\begin{co}
For each pair of integers $r,s>2$, at least one of which is not equal to $4$ and for every pair of $a,b\in\N$,
the double sequence
$$
     \{a T_r(m)+b T_s(n)\}_{m,n\in\N}
$$
misses infinitely many primes 
from every primitive arithmetic progression.
\end{co}

In other words, 
one can only describe the
 primes of the form 
 $am^2+bn^2$
by a finite 
set of congruences
for certain $a,b\in\N$ (see~\cite{C}). However, for 
any pair of $a,b\in\N$, if we replace either $m^2$ or $n^2$ above by another polygonal number, 
the resulting linear combination of polygonal numbers
fails to produce infinitely many primes contained in any given primitive arithmetic progression. 

Theorem~\ref{thm_2}
along with Chebyshev's theorem (Corollary 2.6 in \cite{MV}) and Green-Tao Theorem \cite{GT} implies:

\begin{co}
For each pair of coprime positive integers $a$ and $b$, there are arbitrariliy long arithmetic progressions among the primes of the form $$am^2+bn^2.$$ 
\end{co}
 
It also follows from 
Theorem~\ref{thm_2} and Chebyshev's theorem that the sequence of the primes, which are linear combinations (with coprime positive
integer coefficients) of two polygonal numbers at least one of which is not a square,
has zero relative density with respect to the sequence of all primes. 
So, possible existence of arbitrary long arithmetic progressions among such primes (see the following table)
is undetectable by 
Green-Tao Theorem \cite{GT}.
\
 
\begin{table}[h]
\caption{Examples of arithmetic progressions among the primes which are sums of two triangular numbers}
\label{T1}
 \begin{tabular}{| c | c |} 
\hline
\text{Length}&\makecell{\text{Arithmetic progressions}}\\
\hline
3&31, 37, 43\\
\hline
4&43, 61, 79, 97\\
\hline
5&61, 151, 241, 331, 421\\
\hline
6&7, 37, 67, 97, 127, 157\\
\hline
7&157, 37747, 75337, 112927, 150517, 188107, 225697\\
\hline
8&3307, 69457, 135607, 201757, 267907, 334057, 400207, 466357\\
\hline
9&3823, 6133, 8443, 10753, 13063, 15373, 17683, 19993, 22303\\
\hline
10&\makecell{73091, 135461, 197831, 260201, 322571, 384941, 447311, \\
509681, 572051, 634421}\\
\hline
\end{tabular}
\end{table}

The following table suggests the existence of arbitrarily long arithmetic progressions among the primes which are sums of two triangular number. More generally, we 
ask whether there are arbitrarily long arithmetic progressions among the primes 
which are linear combinations of two polygonal numbers with coprime positive
integer coefficients.
\begin{oq}
For each choice of $a,b,r,s\in\N$ with $r,s>2$ and $\gcd(a,b)=1$, are
there are arbitrarily long arithmetic progressions among the primes of the form 
$
     a T_r(m)+b T_s(n)?
$
\end{oq}

\section{Proofs of the theorems}
The main ingredients of the proof of Theo\textbf{}rem~\ref{thm_2} are a couple of results by
Iwaniec \cite{I}, which we state below.
{\color{white}{\hypertarget{A}{A}}}
\begin{thmA}[Theorem\hspace{2pt}1\hspace{2pt}in\hspace{2pt}\cite{I}]
Let $P(x,y)=ax^2+bxy+cy^2+ex+fy+g$ be an irreducible polynomial
    in $\mathbb{Z}[x,y]$ such that
    $\frac{\partial{P}}{\partial x} $ and $\frac{\partial{P}}{\partial y}$ are linearly independent and $P(x,y)$ represents arbitrarily large odd numbers for \mbox{$x,y\in\N$.}
    Then 
    \begin{align}
        \Bigg(\sum_{\substack{\text{prime\hspace*{3pt}}p\hspace{1pt}\leq N,\\p\hspace{1pt}=P(x,y) \text{\hspace*{3pt}for $x,y\in\N$}}}1\Bigg) \; \; \gg \; \; \frac{N}{\log N},
    \label{2.1}\end{align} if $\Delta:=b^2-4ac$ is a perfect square or $D:=af^2-bef+ce^2+g\Delta=0$,
    and 
    \begin{align}
        \Bigg(\sum_{\substack{\text{prime\hspace*{3pt}}p\hspace{1pt}\leq N,\\p\hspace{1pt}=P(x,y) \text{\hspace*{3pt}for $x,y\in\N$}}}1\Bigg) \; \; \asymp \; \; \frac{N}{(\log N)^{3/2}},
    \label{2.2}\end{align} if $\Delta$ is not a perfect square and $D\neq 0$.
\end{thmA}

\begin{lmB}[Lemma~1 in \cite{I}]
{\color{white}{\hypertarget{B}{B}}}
    Let $P(x,y) \; = \; ax^2+bxy+cy^2+ex+fy+g$.
    Then $\frac{\partial{P}}{\partial x} $ and $\frac{\partial{P}}{\partial y}$ are linearly dependent if and only if $$\Delta \; = \; b^2-4ac \; = \; 0,\ \alpha \; := \; bf-2ce \; = \; 0 \text{ \ and \ } \beta \; := \; be-2af \; = \; 0.$$ 
\end{lmB}

\

\noindent
\emph{Proof of Theorem}~\ref{thm_2}.
For a pair of integers $r,s\geq3$ and for a pair of mutually coprime integers $\mu,\nu>0$, 
considering without loss of generality all possible parities of
the natural numbers $m$ and $n$ in 
\begin{equation}\mu T_r(m)+\nu T_s(n),\nonumber
\end{equation}
we note that (see~\eqref{defpolynum})
the above linear combination 
is equal to 
\begin{align}
    P_1(x,y) \; := \;  2\mu(r-2) x^2+2\nu(s-2) y^2+\mu(4-r) x+\nu (4-s)y
\label{p1}
\end{align}
for $m \; = \; 2x$ and $n \; = \; 2y$,
\begin{align}
    P_2(x,y) \; := \;  2\mu(r-2) x^2+2\nu(s-2) y^2+\mu(4-r) x+\nu sy+\nu
\label{p2}\end{align}
for $m \; = \; 2x$ and $n \; = \; 2y+1$,
\begin{align}
    P_3(x,y) \; := \;  2\mu(r-2) x^2+2\nu(s-2) y^2+\mu rx+\nu sy+\mu+\nu
\label{p3}\end{align}
for $m \; = \; 2x+1$ and $n \; = \; 2y+1$.

\

Since 
the coefficients of $x^2$ and $y^2$
in each of the above polynomials are positive
and since the coefficient 
of $xy$ 
is zero, it follows that 
\begin{equation}\Delta \; < \; 0
\label{negdelta}\end{equation} 
for each of the above polynomials
(see~Lemma~\hyperlink{B}{B}).
Hence, 
$\frac{\partial{P_i}}{\partial x} $ and $\frac{\partial{P_i}}{\partial y}$
linearly independent for all $i\in\{1,2,3\}$.

\

Also, from \eqref{p1}, \eqref{p2} and \eqref{p3}, we note that for each of the polynomials $P_1, P_2$ and $P_3$, we have 
\begin{equation}
       D \; = \; 2\mu\nu^2(r-2)(4-s)^2+2\mu^2\nu(s-2)(4-r)^2,
\label{D}\end{equation}
which equals zero 
if and only if $r=s=4$ (see~Theorem~\hyperlink{A}{A}).

\

In order to
 apply Theorem~\hyperlink{A}{A}, 
we must also 
check 
whether the 
polynomials $P_1,P_2$ $P_3$ are irreducible in 
$\Z[x,y]$.
Suppose, for some $i\in\{1,2,3\}$, there exists $a,b,c,a',b',c'\in\Z$ such that
\begin{align*}
P_i(x,y) \; = \; 
(ax+by+c)(a'x+b'y+c').
\end{align*}
Comparing the coefficients of $x^2$, $y^2$ and $xy$ in both of the above expressions for $P_i(x,y)$, we obtain 
\begin{align*}
    aa' \; = \; 2\mu(r-2),\ bb' \; = \; 2\nu(s-2)\text{ and } ab'+a'b \; = \; 0.
\end{align*}
Since $\mu,\nu>0$ and $r,s>2$, the first two equalities imply that $aa'>0$ and $bb'>0$ and
the last equality implies that there exists $k\in\Q$ such that
$\frac{a}{b}=-\frac{a'}{b'}=k$. Hence,
substituting $kb$ for $a$ and $-kb'$ for $a'$ in the first equality and 
dividing both sides of it
with the corresponding sides
of the second equality above, we obtain
\begin{align*}
    k^2 \; = \; -\frac{\mu(r-2)}{\nu(s-2)}<0,
\end{align*}
 which is absurd!
 Thus, $P_i(x,y)$ is irreducible in $\mathbb{Q}[x,y]$ for all $i\in\{1,2,3\}$. 
Hence, it follows from Gauss's Lemma (see Ch.VI, \S2 in \cite{L}) that $P_1,P_2$ and $P_3$ are
irreducible in $\Z[x,y]$ if and only if these are primitive, i.e. the greatest common
divisor of the coefficients of each of these polynomials is 1. 

\

Since for all $i\in\{1,2,3\}$,
the coefficients of $x^2$ and $y^2$
in $P_i(x,y)$ 
are positive and even, whereas
the coefficient of $xy$ is zero, 
it follows that 
if any of the remaining coefficients of $P_i$ are odd,
then $P_i(x,y)$ 
produces infinitely many
odd numbers for $x,y\in\N$. 
In particular, 
if the greatest common divisors of the coefficients of 
$P_1$, $P_2$,
$P_3$ are 1, then each of these polynomials 
represents arbitrarily large odd numbers. 
In the following 
we compute the gcd $G$ of the coefficients of
$P_1$, $P_2$ and
$P_3$. 
For
$i\in\{1,2,3\}$ and $N\in\N$, 
let
$$\pi_{i}(N) \; := \; \sum_{{\substack{\text{prime\hspace*{3pt}}p\hspace{1pt}\leq N,\\p\hspace{1pt}=P_i(x,y) \text{\hspace*{3pt}for $x,y\in\N$}}}}1.$$

%

\

\noindent\textbf{{Case I.}} For $P_1$, from \eqref{p1} we note that 
$$G \; = \;  \gcd(2\mu(r-2), 2\nu(s-2), \mu(4-r), \nu(4-s)) \; = \; \gcd(4,\mu r,\nu s).$$
Hence, $G=1$ if and only if either $\mu r$ or $\nu s$ is odd.
In particular, $G=1$ implies $D>0$ (see~\eqref{D}). 
So, it follows 
from Theorem~\hyperlink{A}{A} that 
 $$
 \pi_1(N) \; \asymp \; 
        \begin{cases}\dfrac{N}{(\log N)^{3/2}}
         &\text{if $\{\mu r,\nu s\}\nsubseteq2\Z$}\\ \\
1 &\text{otherwise. }
\end{cases}$$

\

\noindent\textbf{{Case II.}} For $P_2$, from \eqref{p2} we note that 
$$G \; = \; \gcd(2\mu(r-2), 2\nu(s-2), \mu(4-r), \nu s, \nu) \; = \; \gcd(4,r,\nu).$$
Hence, $G=1$ if and only if either $r$ or $\nu$ is odd.
So, it follows 
from \eqref{D}, 
Theorem~\hyperlink{A}{A} and Chebyshev's theorem (Corollary 2.6 in \cite{MV}) that 
$$ 
\pi_2(N) \; \asymp \; \begin{cases}
\dfrac{N}{\log N} &\text{if } r=s=4 \text{ and $\nu\notin2\Z$}\\ \\
\dfrac{N}{(\log N)^{3/2}} &\text{if $\{r,s\}\neq\{4\}$ and $\{r,\nu\}\nsubseteq2\Z$}\\ \\
1 &\text{otherwise. }
\end{cases}$$

\

\noindent\textbf{{Case III.}} For $P_3$, from \eqref{p3} we note that 
$$G \; = \; \gcd(2\mu(r-2), 2\nu(s-2), \mu r, \nu s, \mu+\nu) \; = \; \gcd(4,r,s,\mu+\nu).$$
Hence, $G=1$ if and only if either $r$ or $s$ of $\mu+\nu$ is odd.
So, it follows 
from \eqref{D}, 
Theorem~\hyperlink{A}{A} and Chebyshev's theorem (Corollary 2.6 in \cite{MV}) that 
$$ 
\pi_3(N) \; \asymp \; \begin{cases}
\dfrac{N}{\log N} &\text{if } r=s=4 \text{ and $\mu+\nu\notin2\Z$}\\ \\
\dfrac{N}{(\log N)^{3/2}} &\text{if $\{r,s\}\neq\{4\}$ and $\{r,s,\mu+\nu\}\nsubseteq2\Z$}\\ \\
1 &\text{otherwise. }
\end{cases}$$

\

\noindent
Since either $\nu$ or $\mu+\nu$ is odd
for all mutually coprime pairs 
$\mu$ and $\nu$, the claim follows from the above cases. 
\qed

\

\noindent
\emph{Proof of Theorem}~\ref{thm_3}:
Follows immediately from
the proof of the previous theorem 
by considering $P_1(4x,4y\pm1)$ modulo~$4$ or $P_1(4x\pm1,4y)$ modulo~$4$ for $x,y\in\N$,
where $P_1$ is the polynomial defined in $\eqref{p1}$.
\qed

\

\section*{Acknowledgments}

The authors would like to thank Prof. Don Zagier for his comments and suggestions.
All the computations for Table~\ref{T1} was done using PARI/GP. Also, Bhattacharya thankfully acknowledges the support from the Department of Science and Technology, Govt.~of India through the INSPIRE Facul\-ty Award IFA17-MA104 and Rahaman thankfully acknowledges the support from the Ministry of Education, Govt. of India through the Prime Minister Research Fellowship 0501972.


\end{document}